# Sensitivity Based Thevenin Index for Voltage Stability Assessment Considering $N$-1 Contingency

Xiaohu Zhang*, Di Shi*, Xiao Lu†, Zhehan Yi*, Qibing Zhang†, Zhiwei Wang*
*GEIRI North America
†State Grid Jiangsu Electric Power Company
xiaohu.zhang@geirina.net

*Abstract*—This paper proposes an approach to address the voltage stability assessment (VSA) considering $N - 1$ contingency. The approach leverages the sensitivity based Thevenin index (STI) which involves evaluating the Jacobian matrix at current operating condition. Since the $N - 1$ contingency case is hypothetical, there is no information regarding the operating condition after a foreseen contingency. The proposed approach first estimates the post-contingency operating point as well as possible PV-PQ transitions based on the current operating point. Then the STI for each contingency can be predicted using the estimated operating condition. Numerical results based on IEEE 14-bus system demonstrate the accuracy of the proposed approach in predicting the voltage stability margin under contingency. Moreover, the on-line implementation of the proposed approach is promising since it only involves solving several linear equations.

*Index Terms*—Voltage stability assessment, Thevenin equivalence, voltage stability index, $N - 1$ contingency.

## Nomenclature

| | |
|---|---|
| $i, j$ | Index of buses. |
| $\Delta V, \Delta \theta$ | Change of voltage magnitude and angle in two subsequent PMU measurements. |
| $Q^{\min}, Q^{\max}$ | Minimum and maximum reactive power limits for a generator. |
| $V_i, \theta_i$ | Voltage magnitude and angle at bus $i$. |
| $\boldsymbol{f}(\boldsymbol{\theta}, \boldsymbol{V})$ | Active power injection expressions. |
| $\boldsymbol{g}(\boldsymbol{\theta}, \boldsymbol{V})$ | Reactive power injection expressions. |

Other symbols are defined as required in the text. Complex numbers are indicated by an overbar. Matrices and vectors are indicated by bold.

## I. Introduction

DUE to the power market deregulation, the increasing demand of electricity consumptions and the massive integration of renewable energy resources, the aging power grid is under stress [1]. In today's competitive power market, the transmission facilities are often operated close to their security limits, which results in compromised reliability [2]. The construction of new transmission infrastructure is one approach to relieve the transmission burden, but the environmental issues and long construction time make this option difficult [3], [4]. Among the security issues, voltage stability is one of the major concerns for the system operator. It usually starts from a local bus or area but may propagate to a system-wide stability problems [5]. Thus, to maintain a reliable operation of electric power system, it is of great importance for the system operators to accurately and timely assess the voltage stability margin, i.e., the distance between the current operating point and the voltage collapse point.

The model-based approaches, i.e., continuation power flow (CPF) [6] or time-domain simulation [7], are one group of methods to monitor the voltage stability margin. One advantage of the model-based approaches is that they can not only provide the voltage stability margin for the current operating point, but also analyze the impacts brought by the *what-if* scenarios such as $N - 1$ contingencies. Nevertheless, the heavy computational burdens impede the model-based voltage stability assessment (VSA) approaches from on-line applications [2].

In the past decade, the wide deployment of GPS-synchronized phasor measurement units (PMUs) shifts the interests of power industry to the measurement-based VSA. Thevenin equivalence (TE) method is one major approach for the measurement-based VSA. The core idea of TE method is to identify the TE based on measurement data provided by PMUs, i.e., $\bar{V}_i, \bar{I}_i$. With the identified TE, different voltage stability index (VSI) can be derived [8]–[11]. Hence, the calculation of Thevenin equivalence is crucial in the VSA process. A comparative study regarding four approaches to estimate the Thevenin equivalence is provided in [12]. With respect to load area, the authors in [13] merge all the boundary buses in the load area into one fictitious bus and leverage the standard TE method to monitor the voltage stability. In [14], [15], the authors extend the two-bus TE method to three-bus and n+1 bus equivalent system to consider the voltage stability within a load area. Besides of providing the voltage stability margin, the method could also compute the power transfer limit on a specified tie line.

Few of the aforementioned work has addressed the VSA for $N - 1$ transmission contingency in real-time. However, the accurate estimation of the voltage stability margin for the foreseen $N - 1$ contingencies, which is potentially implemented on-line, could enhance the situational awareness of the system operators. A hybrid approach to address the VSA for a series of $N - 1$ contingencies is proposed in [16].

This work is funded by SGCC Science and Technology Program under contract no. 5455HJ160007.

The post-contingency status is first estimated by leveraging the sensitivity method. Then the estimated operating states are treated as fictitious measurements from PMUs and the TE method is applied to achieve the voltage stability margin.

Recently, the authors in [17], [18] propose a sensitivity based Thevenin index (STI) to monitor the static long term voltage stability. The STI is based on wide area measurements and can be leveraged at the control center to validate the local Thevenin index (LTI) calculated by using measurements from a PMU. In this paper, we extend the STI in [17] to address the VSA under a series of foreseen $N-1$ transmission contingencies. The proposed approach only involves solving several linear equations, which makes it promising for on-line implementation.

The remaining sections of this paper are organized as follows. In Section II, the derivation of STI is presented. Section III illustrates details of the proposed VSA method. In Section IV, the numerical results based on IEEE 14-bus system are provided. Finally, some conclusions are given in Section V.

## II. SENSITIVITY BASED THEVENIN INDEX

As mentioned in the introduction section, the STI is proposed in [17]. The derivation of STI is illustrated in this section for completeness.

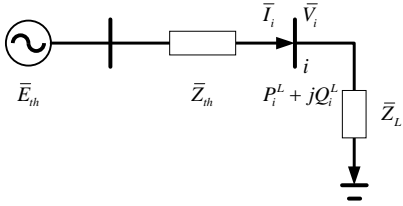

Fig. 1. Two-bus Thevenin equivalent system.

In Fig. 1, a two-bus system with a load bus $i$ and the rest of the system represented by a Thevenin equivalence is provided. To estimate the Thevenin equivalence, i.e., $\bar{E}_{th}, \bar{Z}_{th}$, at least two subsequent measurements are required with the assumption that the Thevenin equivalence remains the same during the time interval. Suppose we have two measurement sets at load bus $i$, i.e., $(\bar{V}_i^{(1)}, \bar{I}_i^{(1)}, \bar{S}_i^{L^{(1)}})$ and $(\bar{V}_i^{(2)}, \bar{I}_i^{(2)}, \bar{S}_i^{L^{(1)}}(1+\Delta\lambda))$, the following two equations hold:

$$\bar{E}_{th} = \bar{V}_i^{(1)} + \bar{Z}_{th} \cdot \bar{I}_i^{(1)} \quad (1)$$
$$\bar{E}_{th} = \bar{V}_i^{(2)} + \bar{Z}_{th} \cdot \bar{I}_i^{(2)} \quad (2)$$

Thus, $\bar{Z}_{th}$ can be calculated as:

$$\bar{Z}_{th} = -\frac{\bar{V}_i^{(2)} - \bar{V}_i^{(1)}}{\bar{I}_i^{(2)} - \bar{I}_i^{(1)}} = -\frac{\Delta \bar{V}_i}{\Delta \bar{I}_i} \quad (3)$$

Based on the voltage stability index proposed in [8], the LTI can be expressed as the ratio between $\bar{Z}_{th}$ and the load impedance $\bar{Z}_L$:

$$\text{LTI} = \left|\frac{\bar{Z}_{th}}{\bar{Z}_L}\right| = \left|\frac{\Delta \bar{V}_i}{\bar{V}_i^{(1)}} \cdot \frac{\bar{I}_i^{(1)}}{\Delta \bar{I}_i}\right| \quad (4)$$

By leveraging the polar form of the complex value and let $\bar{V}_i^{(1)} = Ve^{j\theta}, \bar{V}_i^{(2)} = (V+\Delta V)e^{j(\theta+\Delta\theta)}$, we can write $\Delta \bar{V}_i$ into:

$$\Delta \bar{V}_i = (V+\Delta V)e^{j(\theta+\Delta\theta)} - Ve^{j\theta}$$
$$= \bar{V}_i^{(1)}((1+\Delta V/V)e^{j\Delta\theta} - 1) \quad (5)$$

In addition, $\Delta \bar{I}_i$ can be described by

$$\Delta \bar{I}_i = \bar{I}_i^{(2)} - \bar{I}_i^{(1)} = \frac{\bar{S}_i^{L^{(1)*}}(1+\Delta\lambda)}{Ve^{-j(\theta+\Delta\theta)}} - \frac{\bar{S}_i^{L^{(1)*}}}{Ve^{(-j\theta)}}$$
$$= \bar{I}_i^{(1)}\left(\frac{1+\Delta\lambda}{(1+\Delta V/V)e^{-j\Delta\theta}} - 1\right) \quad (6)$$

With (5) and (6), the ratio between $\bar{Z}_{th}$ and $\bar{Z}_L$ can be derived as:

$$\frac{\bar{Z}_{th}}{\bar{Z}_L} = -\frac{(1+\Delta V/V) \cdot ((1+\Delta V/V) - e^{-j\Delta\theta})}{1+\Delta\lambda - (1+\Delta V/V)e^{-j\Delta\theta}} \quad (7)$$

Since $\Delta\theta$ is small, $e^{-j\Delta\theta} \approx 1 - j\Delta\theta$, the following equation can be obtained:

$$\frac{\bar{Z}_{th}}{\bar{Z}_L} = -(1+\frac{\Delta V}{V}) \cdot \frac{((1+\Delta V/V) - (1-j\Delta\theta))}{1+\Delta\lambda - (1+\Delta V/V)(1-j\Delta\theta)} \quad (8)$$

In (8), $\Delta V$ is small compared to $V$, $\frac{\Delta V}{V}$ can be neglected in the first term. Similarly, the product of $\Delta V$ and $\Delta\theta$ is also small and can be dropped from the denominator in the second term. The simplified version of (8) can be described as:

$$\frac{\bar{Z}_{th}}{\bar{Z}_L} = -\frac{(\Delta V/V + j\Delta\theta)}{\Delta\lambda - \Delta V/V + j\Delta\theta} = -\frac{(\frac{\Delta V}{\Delta\lambda} \cdot \frac{1}{V} + j\frac{\Delta\theta}{\Delta\lambda})}{1 - \frac{\Delta V}{\Delta\lambda} \cdot \frac{1}{V} + j\frac{\Delta\theta}{\Delta\lambda}} \quad (9)$$

Thus, the LTI can be formulated as a function of $\Delta\lambda$:

$$\text{LTI}(\Delta\lambda) = \left|\frac{\bar{Z}_{th}}{\bar{Z}_L}\right| = \sqrt{\frac{(\frac{\Delta V}{\Delta\lambda} \cdot \frac{1}{V})^2 + (\frac{\Delta\theta}{\Delta\lambda})^2}{(1-\frac{\Delta V}{\Delta\lambda} \cdot \frac{1}{V})^2 + (\frac{\Delta\theta}{\Delta\lambda})^2}} \quad (10)$$

The assumption that the Thevenin equivalence remains unchanged during the time interval is only valid when the load increment $\Delta\lambda$ is very small. The idea LTI, i.e., STI, can be achieved by using the limit of (10) with $\Delta\lambda \to 0$:

$$\text{STI} = \lim_{\Delta\lambda \to 0} \text{LTI}(\Delta\lambda) = \sqrt{\frac{(\frac{dV}{d\lambda} \cdot \frac{1}{V})^2 + (\frac{d\theta}{d\lambda})^2}{(1-\frac{dV}{d\lambda} \cdot \frac{1}{V})^2 + (\frac{d\theta}{d\lambda})^2}} \quad (11)$$

To obtain the STI, two sensitivity terms $dV/d\lambda$ and $d\theta/d\lambda$ should be calculated. The procedure can be interpreted as one predictor step of CPF when $d\lambda$ is selected to be the continuation parameter [6], which is provided in equation (12).

$$\begin{bmatrix} f_\theta & f_V & f_\lambda^P \\ g_\theta & g_V & g_\lambda^Q \\ 0 & 0 & 1 \end{bmatrix} \begin{bmatrix} d\theta \\ dV \\ d\lambda \end{bmatrix} = \begin{bmatrix} 0 \\ 0 \\ 1 \end{bmatrix} \quad (12)$$

In (12), $f_\theta, f_V, g_\theta, g_V$ represent the partial derivative of the active and reactive power injection expressions with respect to bus voltage magnitudes and angles, i.e., the standard Jacobian matrix in power flow analysis at a operating point. The term $f_\lambda^P, g_\lambda^P$ indicate the amount of active and reactive power injection changes as a function of $d\lambda$.

## III. PROPOSED APPROACH

Equation (12) should be evaluated at a operating point to obtain the sensitivities. However, the foreseen $N-1$ contingencies are hypothetical and there is no information regarding the operating point. To address this issue, the post-contingency states are estimated based on current operating point. Then equation (12) is evaluated at the estimated operating point to achieve the sensitivities and further the STI for each contingency.

### A. Estimation of Post-contingency Operating States

We leverage the method in [16], [19] to estimate the operating point after a foreseen transmission contingency. The core idea is provided in Fig. 2.

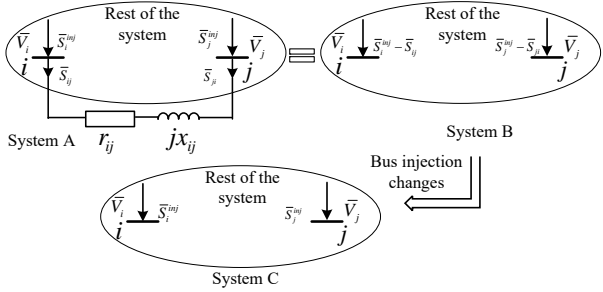

Fig. 2. System transformation.

As shown in Fig. 2, system A represents the power network in base operating condition, where $\bar{S}_i^{inj}$ and $\bar{S}_{ij}$ indicate the net power injection at bus $i$ and line flow from bus $i$ to bus $j$, respectively. Similar description applies to $\bar{S}_j^{inj}$ and $\bar{S}_{ji}$. The difference between system A and system B lies in two parts: 1) system B does not have branch $ij$; 2) the net power injections at bus $i$ and $j$ of system B are changed to be $\bar{S}_i^{inj} - \bar{S}_{ij}$ and $\bar{S}_j^{inj} - \bar{S}_{ji}$, respectively. System A and B are equivalent and have the same power flow results, which is proved in [16]. The post-contingency system, i.e., branch $ij$ in outage, is provided as system C. Compared to system B, the power injections at bus $i$ and $j$ are changed by the amount of $\bar{S}_{ij}$ and $\bar{S}_{ji}$. Thus, the post-contingency operating states can be calculated as:

$$\begin{bmatrix} \boldsymbol{\theta}_c \\ \boldsymbol{V}_c \end{bmatrix} = \begin{bmatrix} \tilde{\boldsymbol{\theta}} \\ \tilde{\boldsymbol{V}} \end{bmatrix} + \begin{bmatrix} \boldsymbol{f_\theta} & \boldsymbol{f_V} \\ \boldsymbol{g_\theta} & \boldsymbol{g_V} \end{bmatrix}^{-1} \Delta \boldsymbol{F_{inj}} \quad (13)$$

where $(\tilde{\boldsymbol{\theta}}, \tilde{\boldsymbol{V}})$ is the operating point under base case; $(\boldsymbol{\theta}_c, \boldsymbol{V}_c)$ is the operating point following a contingency; $\Delta \boldsymbol{F_{inj}}$ is the vector containing the active and reactive bus injection changes which are numerically equal to the active and reactive line flow on the outage branch. Note that the Jacobian matrix in (13) is constructed based on the topology of system B, i.e., the contingency line $ij$ is in outage.

When a severe transmission contingency occurs, some generators may reach their reactive power limits and the corresponding bus type will change from PV to PQ. To both accurately estimate the post-contingency states and the voltage stability margin, it is paramount to predict the possible PV-PQ transitions. The piecewise linear sensitivity method [16] is used to predict the Q limits violation. First, a $K$ factor which reflects the severity of reactive power violation at a PV bus is defined as follows

$$K = \begin{cases} \dfrac{Q^{\max} - \tilde{Q}^{(l)}}{Q^{(l+1)} - \tilde{Q}^{(l)}}, & Q^{(l+1)} \geq Q^{\max} \\ 1, & Q^{\min} \leq Q^{(l+1)} \leq Q^{\max} \\ \dfrac{\tilde{Q}^{(l)} - Q^{\min}}{\tilde{Q}^{(l)} - Q^{(l+1)}}, & Q^{l+1} \leq Q^{\min} \end{cases} \quad (14)$$

where $\tilde{Q}^{(l)}$ and $Q^{(l+1)}$ are the reactive power generation during the iteration process. The complete procedure is provided as below:

1) Set $l = 0$ and obtain the power flow results, $\tilde{\boldsymbol{\theta}}^{(l)}$, $\tilde{\boldsymbol{V}}^{(l)}, \tilde{\boldsymbol{Q}}^{(l)}$, in base operating condition, i.e., before the contingency.
2) Apply $\Delta \boldsymbol{F}_{inj}^{(l)}$ to system $l$ and leverage equation (13) to predict the system $l+1$ states $\boldsymbol{\theta}^{(l+1)}, \boldsymbol{V}^{(l+1)}, \boldsymbol{Q}^{(l+1)}$.
3) Calculate $K$ factor for all the PV buses. If the values of all $K$ are equal to 1 or there is no PV bus, go to step 6); otherwise, go to step 4).
4) Find the PV bus with the lowest $K$, say bus $i$ with $K_i$. Apply a portion of $\Delta \boldsymbol{F}_{inj}^{(l)}$, i.e., $K_i \Delta \boldsymbol{F}_{inj}^{(l)}$ to system $l$ to predict the intermediate system $l+1$ states $\tilde{\boldsymbol{\theta}}^{(l+1)}, \tilde{\boldsymbol{V}}^{(l+1)}, \tilde{\boldsymbol{Q}}^{(l+1)}$.
5) Change bus $i$ to PQ bus by fixing its reactive power generation to its limit. Let $\Delta \boldsymbol{F}_{inj}^{(l+1)} = (1-K_i) \Delta \boldsymbol{F}_{inj}^{(l)}$, set $l = l+1$ and go to step 2).
6) Output the predicted final states $\boldsymbol{\theta}^{(l+1)}, \boldsymbol{V}^{(l+1)}, \boldsymbol{Q}^{(l+1)}$.

### B. Implementation of the Proposed VSA Approach

With the post-contingency states, the STI under each transmission contingency can be predicted. However, equation (12) should be modified to calculate the corresponding sensitivities when there are generators reaching their limits. As an example, suppose bus $i$ is transformed from PV bus to PQ bus. Equation (12) is modified as follows

$$\begin{bmatrix} \boldsymbol{f_\theta} & \boldsymbol{f_V} & \boldsymbol{f}_{V_i} & \boldsymbol{f}_\lambda^P \\ \boldsymbol{g_\theta} & \boldsymbol{g_V} & \boldsymbol{g}_{V_i} & \boldsymbol{g}_\lambda^Q \\ \boldsymbol{g_\theta^i} & \boldsymbol{g_V^i} & g_{V_i}^i & g_\lambda^{Q,i} \\ \boldsymbol{0} & \boldsymbol{0} & 0 & 1 \end{bmatrix} \begin{bmatrix} d\boldsymbol{\theta} \\ d\boldsymbol{V} \\ dV_i \\ d\lambda \end{bmatrix} = \begin{bmatrix} \boldsymbol{0} \\ \boldsymbol{0} \\ 0 \\ 1 \end{bmatrix} \quad (15)$$

Since $V_i$ becomes a variable, one column and row should be added to the previous Jacobian. For the added row, $\boldsymbol{g_\theta^i}, \boldsymbol{g_V^i}, g_{V_i}^i$ represent the partial derivative of reactive power injection expressions at bus $i$ with respect to $\boldsymbol{\theta}, \boldsymbol{V}$ and $V_i$. Note that if bus $i$ does not have any load, $g_\lambda^{Q,i}$ is equal to zero because the reactive limit of the generator at that bus is reached; otherwise, the value of $g_\lambda^{Q,i}$ depends on how will the reactive load varies as $d\lambda$ changes. Equation (15) can be easily extended for any number of PV-PQ bus transitions.

The flowchart of the proposed approach is depicted in Fig. 3. The procedures are described as below:

1) Obtain the current operating states $\tilde{\boldsymbol{\theta}}, \tilde{\boldsymbol{V}}$.

2) Estimate the post-contingency operating point $\boldsymbol{\theta}_c, \boldsymbol{V}_c$ for every considered contingency with the method in Section III-A.
3) If there is no PV-PQ transition following the contingency, calculate the corresponding sensitivities using equation (12); otherwise, obtain the sensitivities using equation (15).
4) Calculate the STI for each considered contingency with equation (11).

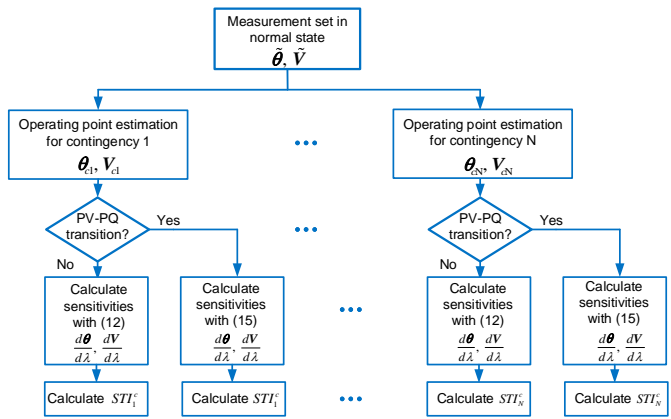

Fig. 3. Flowchart of the proposed VSA approach.

Note that the calculation of STI for each contingency and load bus is independent so the computational speed can be enhanced if parallelization technique is leveraged.

## IV. CASE STUDIES

The IEEE-14 bus system provided in MATPOWER software package [20] is selected to test the performance of the proposed VSA approach. It includes five generators, 14 buses and 20 transmission branches. Since reference [20] only gives one load pattern, we assume that load profile as normal load level, whose active and reactive loads are 259 MW and 73.5 MVAr, respectively. It is also assumed that under peak load level, the active and reactive loads are 1.2 times the values in normal load level. We consider all the $N-1$ branch outages except the following two branches:
- line 1-2 is a key transmission element whose loss will cause the divergence of power flow under normal load pattern;
- line 7-8 only connects a generator bus in one end whose loss will cause islanding system.

First, we implement the proposed approach under normal load condition. Table I provides the VSI, i.e., $\left|\bar{Z}_{th}\right|/\left|\bar{Z}_L\right|$, at every load bus for two critical contingencies. The STI* denotes the benchmark results, which are calculated based on **exact** post-contingency power flow solutions. We compare the accuracy of the proposed approach with the method in [16]. Note that the VSI adopted in [16] is based on the voltage magnitude deviation proposed in [10]. To accommodate the comparison, we shift the VSI to $\left|\bar{Z}_{th}\right|/\left|\bar{Z}_L\right|$. Nevertheless, we leverage the same method as given in [16] to identify the post-contingency TE. The last row in Table I, i.e., $\sigma$, denotes the average relative error (averaged among all load buses) of the two approaches with respect to the benchmark results. As can be observed from Table I, under the normal load level, the proposed approach provides more accurate results than that of the method in [16] for every load bus. The average errors of the proposed method are only 0.31% and 0.96% for contingency 1-5 and 5-6, which are far less than 5.32% and 5.41% given by the method in [16].

TABLE I
COMPARISON OF VSI UNDER NORMAL LOAD LEVEL

| Bus | Contingency Line $(i-j)$ | | | | | |
|---|---|---|---|---|---|---|
| | 1-5 | | | 5-6 | | |
| | STI* | [16] | Proposed | STI* | [16] | Proposed |
| 4 | 0.2874 | 0.2704 | 0.2865 | 0.1976 | 0.1948 | 0.1969 |
| 5 | 0.2780 | 0.2591 | 0.2770 | 0.1570 | 0.1551 | 0.1565 |
| 9 | 0.3590 | 0.3417 | 0.3579 | 0.3577 | 0.3415 | 0.3544 |
| 10 | 0.3638 | 0.3461 | 0.3627 | 0.3763 | 0.3565 | 0.3725 |
| 11 | 0.3695 | 0.3506 | 0.3684 | 0.4085 | 0.3788 | 0.4036 |
| 12 | 0.3821 | 0.3622 | 0.3809 | 0.4443 | 0.4045 | 0.4381 |
| 13 | 0.3797 | 0.3601 | 0.3785 | 0.4401 | 0.4024 | 0.4342 |
| 14 | 0.3784 | 0.3606 | 0.3773 | 0.4092 | 0.3845 | 0.4046 |
| $\sigma$ (%) | - | 5.32 | 0.31 | - | 5.41 | 0.96 |

Fig. 4 depicts the VSIs for all the 18 contingencies at bus 14, which is the critical load bus. The VSI at bus 14 for the base case, i.e., without outage, is 0.2771. It can be seen that the majority of the predicted post-contingency VSI are higher than 0.2771, which indicates that the voltage stability margin will decrease if those contingencies occur. Note that if the VSI exceeds 0.45, alarms will be sent to the control center for situational awareness [17]. Under the normal load level, no alarms will be triggered for bus 14 since the predicted post-contingency VSIs are all below 0.45. Also, the proposed approach outperforms the method in [16] for better estimation of the post-contingency VSI.

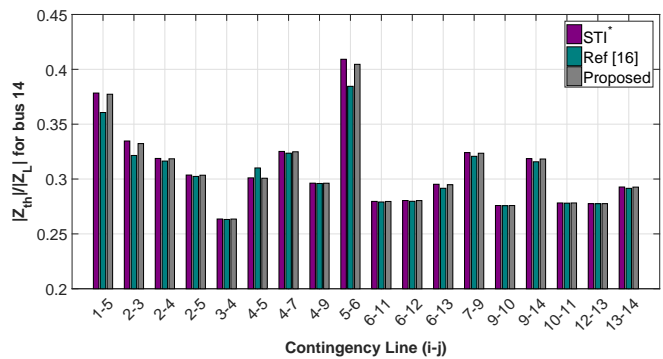

Fig. 4. VSI for different contingencies at bus 14 under normal load level.

Next, we test the proposed approach for the peak load condition. In Table II, the VSI at every load bus for two critical contingencies are provided. Still, the proposed approach gives more accurate results in predicting the post-contingency VSI. For contingency 1-5, the average error $\sigma$ is 24.82% by using the method in [16]. This value decreases to 4.25% by leveraging the proposed approach.

TABLE II
COMPARISON OF VSI UNDER PEAK LOAD LEVEL

| Bus | Contingency Line ($i-j$) | | | | | |
|---|---|---|---|---|---|---|
| | 1-5 | | | 5-6 | | |
| | STI* | [16] | Proposed | STI* | [16] | Proposed |
| 4 | 0.4694 | 0.3345 | 0.4474 | 0.3674 | 0.2081 | 0.3260 |
| 5 | 0.4585 | 0.3111 | 0.4358 | 0.3054 | 0.1642 | 0.2692 |
| 9 | 0.5705 | 0.4389 | 0.5471 | 0.6055 | 0.3229 | 0.5395 |
| 10 | 0.5759 | 0.4441 | 0.5524 | 0.6323 | 0.3477 | 0.5642 |
| 11 | 0.5719 | 0.4377 | 0.5484 | 0.6739 | 0.3784 | 0.6022 |
| 12 | 0.5795 | 0.4441 | 0.5557 | 0.7154 | 0.3555 | 0.6413 |
| 13 | 0.5817 | 0.4473 | 0.5581 | 0.7107 | 0.3566 | 0.6370 |
| 14 | 0.5970 | 0.4666 | 0.5734 | 0.6748 | 0.3669 | 0.6047 |
| $\sigma$ (%) | - | 24.82 | 4.25 | - | 48.25 | 10.82 |

Fig. 5 provides the VSIs for different contingencies at bus 14 under peak load level. As can be seen from Fig. 5, for every contingency, the VSI calculated by the proposed approach is closer to the benchmark results than the VSI given by [16]. The largest difference between the benchmark results and the predicted VSIs occurs in contingency 5-6, where the benchmark VSI is 0.6748 and an alarm should be signaled. The proposed approach calculates the VSI as 0.6047, which still allows the generation of alarm. However, the method in [16] outputs 0.3669 for the VSI, which will not signal any alarm. In addition, the VSI for contingency 1-5, 2-3, 2-4, 5-6 and 7-9 are above 0.45. Thus, the corresponding alarms are generated to trigger some preventive controls, which improves the system operator's situational awareness on the criticality of the current operating condition.

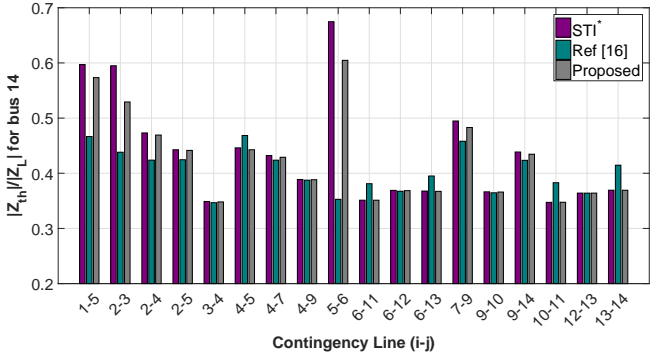

Fig. 5. VSI for different contingencies at bus 14 under peak load level.

## V. CONCLUSIONS

In this paper, an approach based on STI is proposed to address the VSA considering a series of foreseen transmission $N-1$ contingencies. The wide area measurements based STI is interpreted as an idea LTI and can be calculated by evaluating the Jacobian matrix at current operating condition. To extend its application in monitoring the voltage stability margin for $N-1$ contingencies, the post-contingency operating point as well as the PV-PQ transitions are first estimated based on the current operating condition. Then the estimated operating point is used to calculate the STI for each $N-1$ contingency. Numerical case studies based on the IEEE 14-bus system demonstrate the effectiveness of the proposed approach and its benefits in the situational awareness enhancement.